\documentclass[11pt,twoside]{article} 
\date{} 
\title{A note on some infinite sums of Hurwitz zeta functions}
\author{\sc R. B.\ Paris \\
{\em Division of Computing and Mathematics,} \\
{\em Abertay University, Dundee DD1 1HG, UK}}
\begin{document}
\def\f#1#2{\mbox{${\textstyle \frac{#1}{#2}}$}}
\def\dfrac#1#2{\displaystyle{\frac{#1}{#2}}}
\def\boldal{\mbox{\boldmath $\alpha$}}
\newcommand{\bee}{\begin{equation}}
\newcommand{\ee}{\end{equation}}
\newcommand{\sa}{\sigma}
\newcommand{\ka}{\kappa}
\newcommand{\al}{\alpha}
\newcommand{\la}{\lambda}
\newcommand{\ga}{\gamma}
\newcommand{\eps}{\epsilon}
\newcommand{\om}{\omega}
\newcommand{\fr}{\frac{1}{2}}
\newcommand{\fs}{\f{1}{2}}
\newcommand{\g}{\Gamma}
\newcommand{\br}{\biggr}
\newcommand{\bl}{\biggl}
\newcommand{\ra}{\rightarrow}
\newcommand{\gtwid}{\raisebox{-.8ex}{\mbox{$\stackrel{\textstyle >}{\sim}$}}}
\newcommand{\ltwid}{\raisebox{-.8ex}{\mbox{$\stackrel{\textstyle <}{\sim}$}}}
\renewcommand{\topfraction}{0.9}
\renewcommand{\bottomfraction}{0.9}
\renewcommand{\textfraction}{0.05}
\newcommand{\mcol}{\multicolumn}
\date{}
\maketitle
\pagestyle{myheadings}
\markboth{\hfill \sc R. B.\ Paris  \hfill}
{\hfill \sc Hurwitz zeta function sums\hfill}
\begin{abstract}
We consider some closed-form evaluations of certain infinite sums involving the Hurwitz zeta function $\zeta(s,\al)$ of the form
\[\sum_{k=1}^\infty (\pm 1)^k k^m \zeta(s,k),\]
where $m$ is a non-negative integer. For the sums with $m=0$ and the argument $k$ in $\zeta(s,k)$ replaced by $ka+b$, where $a$ and $b$ are positive parameters,  we also obtain a transformation formula suitable for computation in the limit $a\to0$.

\vspace{0.3cm}

\vspace{0.1cm}
 
\end{abstract}

\vspace{0.3cm}

\noindent $\,$\hrulefill $\,$

\vspace{0.3cm}
\begin{center}
{\bf 1.\ Introduction}
\end{center}
\setcounter{section}{1}
\setcounter{equation}{0}
\renewcommand{\theequation}{\arabic{section}.\arabic{equation}}
The Hurwitz zeta function is defined by
\[\hspace{3cm}\zeta(s,\al)=\sum_{n=0}^\infty \frac{1}{(n+\al)^s}\qquad (\Re (s)>1,\ \al\neq 0, -1, -2, \ldots),\]
which reduces to the Riemann zeta function when $a=1$ since $\zeta(s,1)=\zeta(s)$.
We shall make extensive use of the following integral representation for $\zeta(s,\al)$ \cite[(25.11.25)]{DLMF}
\bee\label{e11}
\zeta(s,\al)=\frac{1}{\g(s)} \int_0^\infty \frac{x^{s-1} e^{-\al x}}{1-e^{-x}}dx \qquad (\Re (s)>0,\ \Re (\al)>0)
\ee
together with the large-$\al$ asymptotic behaviour given by \cite[(25.11.43)]{DLMF}
\bee\label{e12}
\zeta(s,\al)= \frac{\al^{1-s}}{s-1} \{1+O(\al^{-1})\}\qquad (\al\to+\infty).
\ee
We shall also employ the binomial expansion expressed in the form
\bee\label{e15}
(1-z)^{-\al}=\sum_{n=0}^\infty \frac{(\al)_n}{n!}\,z^n\qquad (|z|<1),
\ee
where $(\al)_n=\g(\al+n)/\g(\al)$ is the Pochhammer symbol for the rising factorial  and $(1)_n=n!$.
 
Integrals of $\zeta(s,\al)$ with respect to the variable $\al$ in the form
\[\int_0^1f(\al) \zeta(s,\al)\,d\al\]
have been considered in \cite{OM}, for various functions $f(\al)$ including powers of $\al$, exponential, trigonometric and other functions. Integrals involving products of Hurwitz zeta functions, also over the interval $[0,1]$, have been discussed in \cite{SCP}.
  
In Section 9.9 of the monograph by Marichal and Z\'enaidi \cite{MZ}, the authors introduce the function $\ka(s)$  given by
\bee\label{e10}
\ka(s)=\sum_{k=1}^\infty \zeta(s,k)\qquad (\Re (s)>2),
\ee
but claimed that no closed-form evaluation in terms of elementary functions was known. In this note we evaluate this sum, together with its alternating version, in terms of the Riemann zeta function and also consider a variety of other similar sums. A transformation formula is obtained for the sums
\[\sum_{k=0}^\infty (\pm 1)^k \zeta(s,ka+b)\qquad (a>0, b>0)\]
expressed in terms of a similar sum but with the parameter $a$ replaced by $1/a$. This
enables rapid computation in the case $a\to0$, and is analogous to the well-known Poisson-Jacobi transformation for the sum $\sum_{k=0}^\infty \exp\,[-ak^2]$; see \cite[\S 4.1.2]{PK} and Section 4.

\vspace{0.6cm}

\begin{center}
{\bf 2.\ Some closed-form evaluations}
\end{center}
\setcounter{section}{2}
\setcounter{equation}{0}
\renewcommand{\theequation}{\arabic{section}.\arabic{equation}}
We first consider the sum $\ka(s)$ defined in (\ref{e10}) when $\Re (s)>2$. 
Throughout we shall repeatedly carry out inversion of the order of summation and integration, which is permissible by absolute convergence in the stated domain of $\Re (s)$.
From (\ref{e11}), we then have
\[\ka(s)=\frac{1}{\g(s)}\int_0^\infty \frac{x^{s-1}}{1-e^{-x}}\,\sum_{k=1}^\infty e^{-kx}\,dx
=\frac{1}{\g(s)}\int_0^\infty \frac{x^{s-1} e^{-x}}{(1-e^{-x})^2}dx.\]
Inverting the order of summation and integration when $\Re (s)>2$ and employing the binomial expansion of $(1-e^{-x})^{-2}$ in (\ref{e15}), we find
\begin{eqnarray*}
\ka(s)&=&\frac{1}{\g(s)} \sum_{n=0}^\infty \frac{(2)_n}{n!} \int_0^\infty x^{s-1} e^{-(n+1)x} dx=\sum_{n=0}^\infty \frac{(2)_n}{n!}\,(n+1)^{-s}\\
&=&\sum_{n=0}^\infty (n+1)^{1-s},
\end{eqnarray*}
where the integral has been evaluated as a gamma function.
With the change of summation index $n\to n-1$, we therefore obtain the evaluation (see also \cite[(3.10)]{AD})
\bee\label{e21}
\ka(s)=\zeta(s-1)\qquad (\Re (s)>2).
\ee

The alternating version of $\ka(s)$ can be evaluated in a similar manner. We have, when $\Re (s)>1$,
\[{\tilde\ka}(s):=\sum_{k=1}^\infty (-)^{k-1} \zeta(s,k)=\frac{1}{\g(s)}\int_0^\infty \frac{x^{s-1}}{1-e^{-x}}\,\sum_{k=1}^\infty (-)^{k-1}e^{-kx}\,dx,\]
where the sum over $k$ in the integrand equals $e^{-x}/(1+e^{-x})$. Then
\begin{eqnarray*}
{\tilde\ka}(s)&=&\frac{1}{\g(s)}\int_0^\infty \frac{x^{s-1}e^{-x}}{1-e^{-2x}}\,dx=
\frac{1}{\g(s)}\sum_{n=0}^\infty\frac{(1)_n}{n!} \int_0^\infty x^{s-1} e^{-(2n+1)x} dx\\
&=&\sum_{n=0}^\infty (2n+1)^{-s}.
\end{eqnarray*}
Hence we obtain the result
\bee\label{e22}
{\tilde\ka}(s)=(1-2^{-s})\,\zeta(s)\qquad (\Re (s)>1).
\ee

Similar sums that can be evaluated in closed form in terms of Hurwitz zeta functions are
\[\sum_{k=0}^\infty (\pm 1)^k \zeta(s,k+a)\qquad (a>0),\]
where $\Re (s)>2$ (resp. $\Re (s)>1$) for the sum with the upper (resp. lower) sign. For the sum with the upper sign we have
\[\sum_{k=0}^\infty \zeta(s,k+a)=\frac{1}{\g(s)}\int_0^\infty \frac{x^{s-1}}{1-e^{-x}} \sum_{k=0}^\infty e^{-(k+a)x}\,dx=
\frac{1}{\g(s)} \int_0^\infty \frac{x^{s-1} e^{-ax}}{(1-e^{-x})^2}\,dx\]
\[=\frac{1}{\g(s)}\sum_{n=0}^\infty\frac{(2)_n}{n!} \int_0^\infty x^{s-1}e^{-(n+a)x} dx=\sum_{n=0}^\infty (n+1) (n+a)^{-s}.\]
Upon writing $n+1$ as $(n+a)+1-a$, we then obtain
\bee\label{e23}
\sum_{k=0}^\infty \zeta(s,k+a)=\zeta(s-1,a)+(1-a)\,\zeta(s,a)\qquad (\Re (s)>2).
\ee

The alternating version of (\ref{e23}) is
\begin{eqnarray*}
\sum_{k=0}^\infty (-)^k\zeta(s,k+a)&=&\frac{1}{\g(s)}\int_0^\infty \frac{x^{s-1}e^{-ax}}{1-e^{-2x}}\,dx
=\frac{1}{\g(s)}\sum_{n=0}^\infty \frac{(1)_n}{n!}\int_0^\infty x^{s-1} e^{-(2n+a)x}dx\\
&=&\sum_{n=0}^\infty (2n+a)^{-s}.
\end{eqnarray*}
Hence we find that
\bee\label{e24}
\sum_{k=0}^\infty (-)^k\zeta(s,k+a)=2^{-s} \zeta(s, \fs a)\qquad (\Re (s)>1).
\ee
\vspace{0.6cm}

\begin{center}
{\bf 3.\ Moments of Hurwitz zeta function sums}
\end{center}
\setcounter{section}{3}
\setcounter{equation}{0}
\renewcommand{\theequation}{\arabic{section}.\arabic{equation}}
In this section we examine sums of the type
\[\ka(s;m):=\sum_{k=1}^\infty k^m \zeta(s,k) \qquad (m=1, 2, \ldots)\]
subject to the convergence condition $\Re (s)>m+2$, which follows from the asymptotic form (\ref{e12}). The same procedure leads to
the representation
\[\ka(s;m)=\frac{1}{\g(s)} \int_0^\infty \frac{x^{s-1}}{1-e^{-x}}\,\sum_{k=1}^\infty k^m e^{-kx}\,dx.\]
Routine evaluation shows that the inner sum has the form
\[\sum_{k=1}^\infty k^m e^{-kx}=\frac{e^{-x} P_m(e^{-x})}{(1-e^{-x})^{m+1}},\]
where $P_m$ is a polynomial of degree $m-1$. We have for the first few $m$ values:
\[P_1(e^{-x})=1;\quad P_2(e^{-x})=1+e^{-x},\quad P_3(e^{-x})=1+4e^{-x}+e^{-2x},\]
\[P_4(e^{-x})=1+11e^{-x}+11e^{-2x}+e^{-3x}.\]
Then we obtain
\begin{eqnarray*}
\ka(s;m)&=&\frac{1}{\g(s)} \int_0^\infty \frac{x^{s-1} e^{-x}P_m(e^{-x})}{(1-e^{-x})^{m+2}}\,dx\\
&=&\frac{1}{\g(s)}\sum_{n=0}^\infty \frac{(m+2)_n}{n!} \int_0^\infty x^{s-1} e^{-x} P_m(e^{-x})\,dx.
\end{eqnarray*}

We do not evaluate this integral for general integer values of $m$, but simply consider the values $1\leq m\leq 3$.
When $m=1$, we have
\[\ka(s;1)=\frac{1}{\g(s)} \sum_{n=0}^\infty \frac{(3)_n}{n!} \int_0^\infty x^{s-1} e^{-(n+1)x} dx=\frac{1}{2}\sum_{n=0}^\infty (n+2)(n+1)^{1-s}\]
\bee\label{e31}
=\frac{1}{2} \sum_{n=1}^\infty (n+1) n^{1-s}=\frac{1}{2}\{\zeta(s-1)+\zeta(s-2)\}\qquad (\Re (s)>3).
\ee
When $m=2$, we have
\begin{eqnarray}
\ka(s;2)&=&\frac{1}{\g(s)}\sum_{n=0}^\infty \frac{(4)_n}{n!} \int_0^\infty x^{s-1} e^{-(n+1)x} (1+e^{-x})\,dx\nonumber\\
&=&\sum_{n=0}^\infty \frac{(4)_n}{n!}\,\{((n+1)^{-s}+(n+2)^{-s}\}\nonumber\\
&=&\frac{1}{6}\sum_{n=0}^\infty\{(n+2)(n+3) (n+1)^{1-s}+(n+1)(n+3) (n+2)^{1-s}\}\nonumber\\
&=&\frac{1}{6} \sum_{n=1}^\infty \{(n+1)(n+2)+(n^2-1)\} n^{1-s}=\frac{1}{6} \sum_{n=1}^\infty (2n^2+3n+1)n^{1-s}\nonumber\\
&=&\frac{1}{6}\{\zeta(s-1)+3\zeta(s-2)+2\zeta(s-3)\}\qquad (\Re (s)>4).
\end{eqnarray}
Similarly, when $m=3$ we obtain
\begin{eqnarray}
\ka(s;3)&=&\frac{1}{\g(s)}\int_0^\infty \frac{x^{s-1}e^{-x}}{(1-e^{-x})^5}\,(1+4e^{-x}+e^{-2x})\,dx\nonumber\\
&=& \sum_{n=0}^\infty \frac{(5)_n}{n!} \bl\{(n+1)^{-s}+4(n+2)^{-s}+(n+3)^{-s}\br\}\nonumber\\
&=&\frac{1}{24}\sum_{n=1}^\infty (n+1)\bl\{(n+2)(n+3)+4(n-1)(n+2)+(n-1)(n-2)\br\}n^{1-s}\nonumber\\
&=&\frac{1}{4}\sum_{n=1}^\infty (n+1)^2 n^{2-s}\\
&=&\frac{1}{4}\{\zeta(s-2)+2\zeta(s-3)+\zeta(s-4)\}\qquad (\Re (s)>5).
\end{eqnarray}

The alternating variants of these sums are, when $\Re (s)>m+1$,
\[{\tilde\ka}(s;m):=\sum_{k=1}^\infty (-)^{k-1}k^m \zeta(s,k)=\frac{1}{\g(s)}\int_0^\infty \frac{x^{s-1}}{1-e^{-x}}\,\sum_{k=1}^\infty (-)^{k-1}k^m e^{-kx}\,dx,\]
where the inner sums have the form
\[\sum_{k=1}^\infty (-)^{k-1}k^m e^{-kx}=\frac{e^{-x} P_m(-e^{-x})}{(1+e^{-x})^{m+1}}.\]
Thus we obtain the representation
\[{\tilde\ka}(s;m)=\frac{1}{\g(s)}\int_0^\infty \frac{x^{s-1}e^{-x}}{1-e^{-2x}}\,\frac{P_m(-e^{-x})}{(1+e^{-x})^m}\,dx\qquad (\Re (s)>m+1).\]

When $m=1$, we therefore find
\[{\tilde\ka}(s;1)=\frac{1}{\g(s)} \int_0^\infty \frac{x^{s-1} e^{-x}}{(1-e^{-2x})(1+e^{-x})}\,dx
=\frac{1}{\g(s)} \sum_{n=0}^\infty\frac{(1)_n}{n!}\int_0^\infty \frac{x^{s-1} e^{-(2n+1)x}}{1+e^{-x}}\,dx\]
\[=2^{-s}\sum_{n=0}^\infty \{\zeta(s,n+\fs)-\zeta(s,n+1)\}\]
upon use of the result \cite[(25.11.35)]{DLMF}
\bee\label{e300}
\frac{1}{\g(s)} \int_0^\infty \frac{x^{s-1} e^{-\al x}}{1+e^{-x}}\,dx=2^{-s}\{\zeta(s, \fs\al)-\zeta(s,\fs\al+\fs)\}\qquad (\Re (s)>0, \ \Re (\al)>0).
\ee
Then, from (\ref{e21}) and (\ref{e23}) we obtain
\bee\label{e310}
{\tilde\ka}(s;1)=2^{-s}\bl\{\zeta(s-1,\fs)+\fs \zeta(s,\fs)-\zeta(s-1)\br\}\qquad (\Re (s)>2).
\ee

When $m=2$, we have
\begin{eqnarray}
{\tilde\ka}(s;2)&=&\frac{1}{\g(s)}\int_0^\infty \frac{x^{s-1}e^{-x}}{(1+e^{-x})^3}\,dx
=\frac{1}{\g(s)}\sum_{n=0}^\infty\frac{(-)^n (3)_n}{n!} \int_0^\infty x^{s-1} e^{-(n+1)x}\,dx\nonumber\\
&=&\frac{1}{2}\sum_{n=0}^\infty \frac{(-)^n (3)_n}{n!} (n+1)^{-s} =\frac{1}{2}\sum_{n=1}^\infty (-)^{n-1}(n+1)n^{1-s}\nonumber\\
&=&\frac{1}{2}\bl\{(1-2^{2-s}) \zeta(s-1)+(1-2^{3-s}) \zeta(s-2)\br\}\qquad (\Re (s)>3)\label{e32}
\end{eqnarray}
upon use of the standard result \cite[(25.2.3 )]{DLMF}
\[\sum_{n=1}^\infty (-)^{n-1} n^{-s}=(1-2^{1-s}) \zeta(s)\qquad (\Re (s)>0).\] 
We have been unable to evaluate the case $m=3$ in closed form.

The alternating sums can also be written in the form
\[{\tilde\ka}(s;m)=\ka(s;m)-2\sum_{k=1}^\infty (2k)^m \zeta(s,2k),\]
whence it follows that
\[\sum_{k=1}^\infty k^m \zeta(s,2k)=2^{-m-1}\{\ka(s;m)-{\tilde\ka}(s;m)\}.\]
From (\ref{e31}) and (\ref{e32}), we therefore find the evaluations
\[\sum_{k=1}^\infty k\,\zeta(s,2k)=\frac{1}{8}\bl\{(1+2^{1-s})\zeta(s-1)+\zeta(s-2)-2^{1-s}\bl(\zeta(s-1,\fs)+\fs \zeta(s,\fs)\br)\br\}\]
when $\Re (s)>3$ and
\[\sum_{k=1}^\infty k^2\,\zeta(s,2k)=\frac{1}{24}\bl\{(3\cdot 2^{1-s}-1) \zeta(s-1)+6\cdot 2^{1-s} \zeta(s-2)+\zeta(s-3)\br\}\]
when $\Re (s)>4$.
\vspace{0.6cm}

\begin{center}
{\bf 4.\  Two transformation formulas}
\end{center}
\setcounter{section}{4}
\setcounter{equation}{0}
\renewcommand{\theequation}{\arabic{section}.\arabic{equation}}
In this final section we consider the sums
\bee\label{e41}
\sum_{k=0}^\infty (\pm 1)^k \zeta(s,ka+b)\qquad (a>0, b>0),
\ee
where $\Re (s)>2$ (resp. $\Re (s)>1$) for the series with the upper (resp. lower) sign.
Although we do not achieve a closed-form evaluation of these sums, we express them in terms of related sums where the parameter $a$ is inverted. This results in a convenient means of computation in the case $a\to0$, where convergence of (\ref{e41}) will be slow, since the terms of the transformed sums contain $k/a$ and so decay rapidly with increasing $k$. Such a result is analogous to the well-known Poisson-Jacobi transformation for the sum of Gaussian exponentials given by \cite[\S 4.1.2]{PK}, \cite[p.~124]{WW}
\[\sum_{k=1}^\infty e^{-ak^2}=\frac{1}{2}\sqrt{\frac{\pi}{a}}-\frac{1}{2}+\sqrt{\frac{\pi}{a}} \sum_{k=1}^\infty e^{-\pi^2k^2/a}.\]

Consider first (\ref{e41}) with the positive sign. We have
\[\ka(s;a,b):=\sum_{k=0}^\infty \zeta(s,ka+b)=\frac{1}{\g(s)} \int_0^\infty \frac{x^{s-1}}{1-e^{-x}}\sum_{k=0}^\infty e^{-(ka+b)x} dx\]
\[=\frac{1}{\g(s)}\int_0^\infty \frac{x^{s-1}e^{-bx}}{(1-e^{-x})(1-e^{-ax})}\,dx.\]
If we now expand the factor $(1-e^{-x})^{-1}$ by the binomial theorem, we obtain
\[\ka(s;a,b)=\frac{1}{\g(s)}\sum_{n=0}^\infty \int_0^\infty\frac{x^{s-1} e^{-(n+b)x}}{1-e^{-ax}}\,dx
=\frac{a^{-s}}{\g(s)} \sum_{n=0}^\infty \int_0^\infty \frac{t^{s-1} e^{-(n+b)t/a}}{1-e^{-t}}\,dt\]
\bee\label{e42}
=a^{-s} \sum_{n=0}^\infty \zeta\bl(s,\frac{n+b}{a}\br)\qquad (\Re (s)>2).
\ee

The treatment of the alternating version
\[{\tilde\ka}(s;a,b)=\sum_{k=0}^\infty (-)^k \zeta(s,ka+b)\]
follows in a similar manner. We have
\[{\tilde\ka}(s;a,b)=\frac{1}{\g(s)} \int_0^\infty \frac{x^{s-1}}{1-e^{-x}} \sum_{n=0}^\infty (-)^k e^{-(ka+b)x}\,dx
=\frac{1}{\g(s)} \int_0^\infty\frac{x^{s-1} e^{-bx}}{(1-e^{-x})(1+e^{-ax})}\,dx\]
\[=\frac{1}{\g(s)}\sum_{n=0}^\infty\int_0^\infty\frac{x^{s-1}e^{-(n+b)x}}{1+e^{-ax}}\,dx.\]
Making the change of variable $t=ax$ and employing the integral (\ref{e300}), we find
\bee\label{e43}
{\tilde\ka}(s;a,b)=(2a)^{-s} \sum_{n=0}^\infty \bl\{\zeta\bl(s,\frac{n+b}{2a}\br)-\zeta\bl(s,\frac{n+b}{2a}+\frac{1}{2}\br)\br\}\qquad (\Re (s)>1).
\ee
\bigskip

\noindent{\bf Corollary.}\ \ When $b=a$, we obtain from (\ref{e42}) and (\ref{e43}), upon replacing the summation index $k$ by $k-1$,
\[\sum_{k=1}^\infty \zeta(s,ka)=a^{-s}\sum_{n=0}^\infty \zeta\bl(s,\frac{n}{a}+1\br)\qquad(\Re (s)>2),\]
\[\sum_{k=1}^\infty (-)^{k-1} \zeta(s,ka)=(2a)^{-s}\sum_{n=0}^\infty\bl\{\zeta\bl(s,\frac{n}{2a}+\frac{1}{2}\br)-\zeta\bl(s,\frac{n}{2a}+1\br)\br\}\qquad(\Re (s)>1).\]

As a numerical example, consider the computation of $\ka(s;a,b)$ when $s=4$, $b=1$. The terms in the sum decay like $(ka)^{-3}$ with the result that computation of the terms down to $O(10^{-8})$, say,  when $a=10^{-1}$ requires approximately $1.5\times 10^3$ terms. The number of terms required on the right-hand side of (\ref{e42}) to the same accuracy is about 14. When $a=10^{-2}$, the situation is more acute, requiring about $1.5\times 10^4$ terms for the left-hand side of (\ref{e42}), but only 2 terms for the right-hand side. 

Finally, the sums
\[S_\pm(s;a,b,c)=\sum_{k=0}^\infty (\pm 1)^k e^{-ck} \zeta(s,ka+b)\qquad (\Re (s)>1;\  a, b, c>0)\]
can also be expressed in a similar form with an inversion of the parameter $a$. We have from (\ref{e11})
\[S_\pm(s;a,b,c)=\frac{1}{\g(s)} \int_0^\infty \frac{x^{s-1}}{1-e^{-x}}\,\sum_{k=0}^\infty (\pm 1)^k e^{-(ka+b)x-ck}dx\]
\[\hspace{1cm}=\frac{1}{\g(s)}\int_0^\infty \frac{x^{s-1} e^{-bx}}{(1-e^{-x})(1\mp e^{-ax-c})}\,dx.\]
Expanding the factor $(1-e^{-x})^{-1}$ by the binomial theorem, we obtain
\[S_\pm(s;a,b,c)=\frac{1}{\g(s)}\sum_{n=0}^\infty \frac{(1)_n}{n!} \int_0^\infty \frac{x^{s-1} e^{-(n+b)x}}{1\mp e^{-ax-c}}\,dx
=\frac{a^{-s}}{\g(s)} \sum_{n=0}^\infty\int_0^\infty \frac{t^{s-1} e^{-(n+b)t/a}}{1\mp e^{-t-c}}\,dx.\]

The integral can be evaluated as a Lerch transcendent $\Phi(z,s,\al)$ defined by \cite[p.~612]{DLMF} 
\[\hspace{2cm}\Phi(z,s,\al)=\sum_{n=0}^\infty \frac{z^n}{(n+\al)^s}\qquad (|z|<1; \ \Re(s)>1, |z|=1)\]
when $\al\neq 0, -1, -2, \ldots$, and possessing the integral representation
\[\hspace{2cm}\Phi(z,s,\al)=\int_0^\infty \frac{x^{s-1} e^{-\al x}}{1-z e^{-x}}\,dx\qquad (\Re (s)>0, \, \Re (\al)>0,\ z\in {\bf C}\backslash [1,\infty)).\]
Thus, upon identifying $z$ as $\pm e^{-c}$, we obtain the transformation formula
\bee\label{e44}
S_\pm(s;a,b,c)=a^{-s} \sum_{n=0}^\infty \Phi\bl(\pm e^{-c},s,\frac{n+b}{a}\br).
\ee
When $c=0$ the right-hand sides of (\ref{e44}) reduce to (\ref{e42}) and (\ref{e43}) since $\Phi(1,s,\al)=\zeta(s,\al)$ and $\Phi(-1,s,\al)=2^{-s}\{\zeta(s,\fs \al)-\zeta(s,\fs\al+\fs)\}$ when $\Re (s)>1$.

\end{document}